# SAMPLED-DATA BOUNDARY FEEDBACK CONTROL OF 1-D PARABOLIC PDES

**Iasson Karafyllis[*] and Miroslav Krstic[**]**

[*]Dept. of Mathematics, National Technical University of Athens,
Zografou Campus, 15780, Athens, Greece, email: iasonkar@central.ntua.gr

[**]Dept. of Mechanical and Aerospace Eng., University of California, San Diego, La Jolla, CA 92093-0411, U.S.A., email: krstic@ucsd.edu

**Abstract**

The paper provides results for the application of boundary feedback control with Zero-Order-Hold (ZOH) to 1-D linear parabolic systems on bounded domains. It is shown that the continuous-time boundary feedback applied in a sample-and-hold fashion guarantees closed-loop exponential stability, provided that the sampling period is sufficiently small. Two different continuous-time feedback designs are considered: the reduced model design and the backstepping design. The obtained results provide stability estimates for weighted 2-norms of the state and robustness with respect to perturbations of the sampling schedule is guaranteed.

## 1. Introduction

Sampled-data feedback control is a well-studied topic for finite-dimensional systems due to the fact that modern control systems employ digital technology for the implementation of the controller (see for instance [5,9,10,11,22,23,32] and the references therein). However, for infinite-dimensional systems there are few results on sampled-data feedback control. Most of the available results deal with delay systems (see [6,13,25,26,27,28,31]). For systems described by Partial Differential Equations (PDEs) the design of sampled-data feedback control is a major challenge because of the technical difficulties involved: even the notion of the solution of a PDE under sampled-data feedback control has to be clarified. Sampled-data controllers for parabolic systems were designed by Fridman and coworkers in [1,7,8,30] by using matrix inequalities. The major works [19,29] provided necessary and sufficient conditions for sampled-data control of general infinite-dimensional systems under periodic sampling. The results in [19,29] were extended in the case of "generalized sampling" in [20,35]. Approximate models of infinite-dimensional systems were used in [34] for practical stabilization. Sampled-data feedback control for hyperbolic age-structured models was proposed in [14].

In the linear finite-dimensional case, there are results that guarantee closed-loop exponential stability for continuous-time linear feedback designs when applied with Zero-Order-Hold (ZOH) and sufficiently small sampling period (see for instance [9,10,22,23]). The results deal with the globally Lipschitz case (which contains the linear case as a subcase) and the application of the continuous-time feedback under ZOH is called the "emulation" sampled-data feedback design.

A general robustness result that guarantees closed-loop exponential stability for continuous-time linear boundary feedback designs when applied with ZOH and arbitrary (not necessarily periodic) sampling schedules of sufficiently small sampling period is missing for the case of linear parabolic PDEs. The present work is devoted to the development of such general results for continuous-time boundary feedback controllers. We consider two different boundary feedback designs for linear parabolic PDE models:
- the "Reduced Model Design", originally proposed in [36] and later studied in [3,4,17,18], and
- the "Backstepping Design", developed in [16,33].

In each of the above cases, we prove that there is a sufficiently small sampling period, such that the closed-loop system preserves exponential stability under sample-and-hold implementation of the controller (Theorem 3.1 and Theorem 3.2). In other words, we prove that emulation design works for the case of linear parabolic PDEs with boundary feedback. The derived exponential stability estimates are expressed in appropriate weighted $L^2$ norms of the state and (conservative) upper bounds for the sampling period are derived. Finally, robustness with respect to the sampling schedule is established, exactly as in the finite-dimensional case.

The methodology for proving the main results of the present work is very different for each boundary feedback design. The reason for the difference in the proofs is that the reduced model design and the backstepping design are very different (although they look similar; see Section 3.III for additional explanations). Another thing that should be emphasized at this point is that a general Lyapunov-like proof that can work for every linear parabolic PDE under a linear stabilizing boundary feedback is not available and may not exist (contrary to the finite-dimensional case): the effect of the boundary input is expressed by unbounded linear operators. Consequently, the effect of the difference between the continuous-time feedback and the applied control action is also expressed by linear unbounded operators, which makes a Lyapunov-like analysis of the closed-loop system difficult.

The structure of the present work is as follows: Section 2 is devoted to the presentation of the problem and the clarification of the notion of the solution for a parabolic system under boundary sampled-data control. Section 3 provides the statements of the main results (Theorem 3.1 and Theorem 3.2) as well as a discussion of the main results. The proofs of the main results are provided in Section 4. A simple illustrating example is presented in Section 5. The concluding remarks are provided in Section 6. Finally, the Appendix contains the proofs of all auxiliary results.

**Notation.** Throughout this paper, we adopt the following notation.

* $\Re_+ := [0,+\infty)$. $Z^+$ denotes the set of all non-negative integers.
* Let $U \subseteq \Re^n$ be a set with non-empty interior and let $\Omega \subseteq \Re$ be a set. By $C^0(U;\Omega)$, we denote the class of continuous mappings on $U$, which take values in $\Omega$. By $C^k(U;\Omega)$, where $k \geq 1$, we denote the class of continuous functions on $U$, which have continuous derivatives of order $k$ on $U$ and take values in $\Omega$.
* For a vector $x \in \Re^m$, $x'$ denotes the transpose of $x$. For real numbers $a_i$, $i=1,...,m$, $diag(a_1,a_2,...,a_m)$ denotes the diagonal square matrix with $a_1, a_2,...,a_m$ on its main diagonal.
* Let $r \in C^0([0,1];(0,+\infty))$ be given. $L_r^2([0,1])$ denotes the equivalence class of measurable functions $f:[0,1] \to \Re$ for which $\|f\|_r = \left( \int_0^1 r(z)|f(z)|^2 dz \right)^{1/2} < +\infty$. $L_r^2([0,1])$ is a Hilbert space with inner product $\langle f,g \rangle = \int_0^1 r(z)f(z)g(z)dz$. When $r(z) \equiv 1$, we use the notation $L^2([0,1])$ for the standard space of square-integrable functions and $\|f\|_2 = \left( \int_0^1 |f(z)|^2 dz \right)^{1/2} < +\infty$ for $f \in L^2([0,1])$.
* Let $x:\Re_+ \times [0,1] \to \Re$ be given. We use the notation $x[t]$ to denote the profile at certain $t \geq 0$, i.e., $(x[t])(z) = x(t,z)$ for all $z \in [0,1]$.
* $H^2(0,1)$ denotes the Sobolev space of continuously differentiable functions on $[0,1]$ with measurable, square integrable second derivative.

## 2. The notion of the Solution for a Parabolic System Under Boundary Sampled-Data Control

Consider the Sturm-Liouville operator $A: D \to C^0([0,1]; \Re)$ defined by

$$(Af)(z) = -\frac{1}{r(z)} \frac{d}{dz}\left(p(z)\frac{df}{dz}(z)\right) + \frac{q(z)}{r(z)} f(z), \text{ for all } f \in D \text{ and } z \in (0,1) \tag{2.1}$$

where $p \in C^1([0,1]; (0,+\infty))$, $r \in C^0([0,1]; (0,+\infty))$, $q \in C^0([0,1]; \Re)$ and $D \subseteq H^2([0,1])$ is the set of all functions $f:[0,1] \to \Re$ for which

$$b_1 f(0) + b_2 \frac{df}{dz}(0) = a_1 f(1) + a_2 \frac{df}{dz}(1) = 0 \tag{2.2}$$

where $a_1, a_2, b_1, b_2$ are real constants with $|a_1| + |a_2| > 0$, $|b_1| + |b_2| > 0$. The following fact is a direct consequence of Chapter 11 in [2] and pages 498-505 in [21].

**FACT:** All eigenvalues of the Sturm-Liouville operator $A: D \to L^2_r([0,1])$, defined by (2.1), (2.2) are real. They form an infinite, increasing sequence $\lambda_1 < \lambda_2 < \ldots < \lambda_n < \ldots$ with $\lim_{n \to \infty}(\lambda_n) = +\infty$. To each eigenvalue $\lambda_n \in \Re$ ($n = 1, 2, \ldots$) corresponds exactly one eigenfunction $\phi_n \in C^2([0,1]; \Re)$ that satisfies $A\phi_n = \lambda_n \phi_n$ and $b_1 \phi_n(0) + b_2 \frac{d\phi_n}{dz}(0) = a_1 \phi_n(1) + a_2 \frac{d\phi_n}{dz}(1) = 0$. The eigenfunctions form an orthonormal basis of $L^2_r([0,1])$.

In the present work, we make the following assumption for the Sturm-Liouville operator $A: D \to C^0([0,1]; \Re)$ defined by (2.1), (2.2), where $a_1, a_2, b_1, b_2$ are real constants with $|a_1| + |a_2| > 0$, $|b_1| + |b_2| > 0$.

**(H):** *The Sturm-Liouville operator $A: D \to L^2_r([0,1])$ defined by (2.1), (2.2), where $a_1, a_2, b_1, b_2$ are real constants with $|a_1| + |a_2| > 0$, $|b_1| + |b_2| > 0$, satisfies*

$$\sum_{n=N}^{\infty} \lambda_n^{-1} \max_{0 \le z \le 1}(|\phi_n(z)|) < +\infty, \text{ for certain } N > 0 \text{ with } \lambda_N > 0 \tag{2.3}$$

We consider the following control system

$$\frac{\partial x}{\partial t}(t,z) - \frac{1}{r(z)} \frac{\partial}{\partial z}\left(p(z)\frac{\partial x}{\partial z}(t,z)\right) + \frac{q(z)}{r(z)} x(t,z) = 0, \ z \in (0,1) \tag{2.4}$$

$$b_1 x(t,0) + b_2 \frac{\partial x}{\partial z}(t,0) = a_1 x(t,1) + a_2 \frac{\partial x}{\partial z}(t,1) - u(t) = 0, \tag{2.5}$$

where $x[t]$ is the state and $u(t)$ is the boundary control input. More specifically, we consider system (2.4), (2.5) under boundary sampled-data control with ZOH:

$$u(t) = u_i, \text{ for } t \in [\tau_i, \tau_{i+1}) \text{ and for all } i \in Z^+ \tag{2.6}$$

where $\{\tau_i \geq 0, i = 0,1,2,...\}$ is an increasing sequence (the sequence of sampling times) with $\tau_0 = 0$, $\lim_{i \to +\infty}(\tau_i) = +\infty$ and $\{u_i \in \Re, i = 0,1,2,...\}$ is the sequence of applied inputs.

In order to study system (2.4), (2.5) under boundary sampled-data control, we first need to clarify the notion of the solution for this system. To this purpose the following theorem is a crucial result. Its proof is given in the Appendix.

**Theorem 2.1:** *Consider the Sturm-Liouville operator $A: D \to L_r^2([0,1])$ defined by (2.1), (2.2), where $a_1, a_2, b_1, b_2$ are real constants with $|a_1| + |a_2| > 0$, $|b_1| + |b_2| > 0$, under Assumption (H). Let $T > 0$ be a constant and let $u:[0,T) \to \Re$ be a right continuous function which is $C^2$ on $(0,T)$ with all meaningful right and left limits of $u(t), \dot{u}(t), \ddot{u}(t)$ when $t$ tends to $0$ or $T$ being finite. Then for every $x_0 \in L_r^2([0,1])$, there exists a unique function $x: \Re_+ \times [0,1] \to \Re$ for which the mapping $[0,T] \ni t \to x[t] \in L_r^2([0,1])$ is continuous, with $x \in C^1((0,T) \times [0,1])$ satisfying $x[t] \in C^2([0,1])$ for all $t \in (0,T]$, $x(0,z) = x_0(z)$ for all $z \in [0,1]$ and*

$$\frac{\partial x}{\partial t}(t,z) - \frac{1}{r(z)} \frac{\partial}{\partial z}\left(p(z) \frac{\partial x}{\partial z}(t,z)\right) + \frac{q(z)}{r(z)} x(t,z) = 0, \text{ for all } (t,z) \in (0,T) \times (0,1) \quad (2.7)$$

$$b_1 x(t,0) + b_2 \frac{\partial x}{\partial z}(t,0) = a_1 x(t,1) + a_2 \frac{\partial x}{\partial z}(t,1) - u(t) = 0, \text{ for all } t \in (0,T) \quad (2.8)$$

Having clarified how the solution can be constructed in an interval, we are in a position to obtain the solution of system (2.4), (2.5) under boundary sampled-data control. The following corollary is a direct consequence of Theorem 2.1 (applied in a step-by-step process).

**Corollary 2.2:** *Consider the Sturm-Liouville operator $A: D \to L_r^2([0,1])$ defined by (2.1), (2.2), where $a_1, a_2, b_1, b_2$ are real constants with $|a_1| + |a_2| > 0$, $|b_1| + |b_2| > 0$, under Assumption (H). Let $\{\tau_i \geq 0, i = 0,1,2,...\}$ be an increasing sequence of sampling times with $\tau_0 = 0$, $\lim_{i \to +\infty}(\tau_i) = +\infty$ and let $\{u_i \in \Re, i = 0,1,2,...\}$ be a sequence of applied inputs. Then for every $x_0 \in L_r^2([0,1])$ there exists a unique function $x: \Re_+ \times [0,1] \to \Re$ for which the mapping $\Re_+ \ni t \to x[t] \in L_r^2([0,1])$ is continuous, with $x \in C^1(I \times [0,1])$ satisfying $x[t] \in C^2([0,1])$ for all $t > 0$, $x(0,z) = x_0(z)$ for all $z \in [0,1]$, and*

$$\frac{\partial x}{\partial t}(t,z) - \frac{1}{r(z)} \frac{\partial}{\partial z}\left(p(z) \frac{\partial x}{\partial z}(t,z)\right) + \frac{q(z)}{r(z)} x(t,z) = 0, \text{ for all } (t,z) \in I \times (0,1) \quad (2.9)$$

$$b_1 x(t,0) + b_2 \frac{\partial x}{\partial z}(t,0) = a_1 x(t,1) + a_2 \frac{\partial x}{\partial z}(t,1) - u(t) = 0, \text{ for all } t \in I \quad (2.10)$$

*where $I = \Re_+ \setminus \{\tau_i \geq 0, i = 0,1,2,...\}$ and $u(t)$ satisfies (2.6).*

It should be noticed that exactly as in the finite-dimensional case the solution of the system is continuous at each time instant and differentiable only in the open intervals that are formed between two consecutive sampling times. However, the solution mapping is only continuous as a mapping in $L_r^2(0,1)$. The solution will not be (in general) a continuous mapping in $L^\infty(0,1)$ in the case of Dirichlet actuation (i.e., if $a_2 = 0$) and if the sequence of applied inputs is not constant. It is therefore clear that in the infinite-dimensional case, the choice of the state space is crucial for the appropriate notion of the solution.

# 3. Stability Under Sampled-Data Implementation

### 3.I. Emulation of the Reduced Model Design

The reduced model boundary feedback design is a method proposed in [3,4,17,18,36] in which all unstable eigenvalues of the Sturm-Liouville operator $A: D \to L_r^2([0,1])$ defined by (2.1), (2.2) are included in the linear finite-dimensional system

$$\dot{x}_n(t) + \lambda_n x_n(t) = g_n u(t), \quad n = 1, \ldots, m \tag{3.1}$$

where $g_n$ ($n = 1, \ldots, m$) are real constants, $g = (g_1, \ldots, g_m)' \in \Re^m$ and $x(t) = (x_1(t), \ldots, x_m(t))' \in \Re^m$. The integer $m \geq 1$ is selected to be sufficiently large so that $\lambda_{m+1} > 0$. The continuous-time feedback is selected to be $u(t) = \int_0^1 r(z) \left( \sum_{l=1}^m k_l \phi_l(z) \right) x(t,z) dz$, where $k = (k_1, \ldots, k_m)' \in \Re^m$ is a vector for which the matrix $diag(-\lambda_1, \ldots, -\lambda_m) + gk'$ is Hurwitz. The following result shows that the boundary sampled-data feedback law obtained by the emulation of the continuous-time reduced model design works.

**Theorem 3.1:** *Consider the Sturm-Liouville operator $A: D \to L_r^2([0,1])$ defined by (2.1), (2.2), where $a_1, a_2, b_1, b_2$ are constants with $|a_1| + |a_2| > 0$, $|b_1| + |b_2| > 0$, under Assumption (H). Let $m \geq 1$ be an integer so that $\lambda_{m+1} > 0$. Let $k = (k_1, \ldots, k_m)' \in \Re^m$ be a vector for which the matrix $diag(-\lambda_1, \ldots, -\lambda_m) + gk'$ is Hurwitz, where $g = (g_1, \ldots, g_m)' \in \Re^m$ and $g_n := \frac{p(1)}{a_1^2 + a_2^2} \left( a_2 \phi_n(1) - a_1 \frac{d\phi_n}{dz}(1) \right) \neq 0$ for $n = 1, \ldots, m$. Then there exist constants $G, c, T > 0$ such that for every increasing sequence $\{\tau_i \geq 0, i = 0,1,2,\ldots\}$ with $\tau_0 = 0$, $\sup_{i \geq 0}(\tau_{i+1} - \tau_i) \leq T$, $\lim_{i \to +\infty}(\tau_i) = +\infty$ and for every $x_0 \in L_r^2([0,1])$ the unique function $x: \Re_+ \times [0,1] \to \Re$ for which the mapping $\Re_+ \ni t \to x[t] \in L_r^2([0,1])$ is continuous, with $x \in C^1(I \times [0,1])$ satisfying $x[t] \in C^2([0,1])$ for all $t > 0$, $x(0,z) = x_0(z)$ for all $z \in [0,1]$, (2.9), (2.10) with $I = \Re_+ \setminus \{\tau_i \geq 0, i = 0,1,2,\ldots\}$ and*

$$u(t) = k' \begin{bmatrix} \int_0^1 r(z) \phi_1(z) x(\tau_i, z) dz \\ \vdots \\ \int_0^1 r(z) \phi_m(z) x(\tau_i, z) dz \end{bmatrix} = \int_0^1 r(z) \left( \sum_{l=1}^m k_l \phi_l(z) \right) x(\tau_i, z) dz, \text{ for } t \in [\tau_i, \tau_{i+1}) \text{ and for all } i \in Z^+, \tag{3.2}$$

*satisfies the following estimate*

$$\|x[t]\|_r \leq G \exp(-ct) \|x_0\|_r, \text{ for all } t \geq 0. \tag{3.3}$$

### 3.II. Emulation of the Backstepping Design

Here, we focus for simplicity reasons on the following control system

$$\frac{\partial x}{\partial t}(t,z) = p \frac{\partial^2 x}{\partial z^2}(t,z) - q(z) x(t,z), \tag{3.4}$$

$$b_1 x(t,0) + b_2 \frac{\partial x}{\partial z}(t,0) = x(t,1) - u(t) = 0, \tag{3.5}$$

where $q \in C^0([0,1]; \Re)$, $p > 0, b_1, b_2 \in \Re$ are constants with $|b_1| + |b_2| > 0$, $x[t]$ is the state and $u(t)$ is the control input. Indeed, there exists a transformation that brings system (2.4) to the form (3.4) (see page 44 in [16]). The reason for considering only the Dirichlet actuation case will be explained shortly.

Theorem 2 in [33] guarantees for every $c \geq 0$ the existence of functions $K, L \in C^2([0,1]^2)$ such that the Volterra transformation

$$y(t, z) = x(t, z) - \int_0^z K(z, s) x(t, s) ds, \text{ for all } (t, z) \in \Re_+ \times [0,1] \qquad (3.6)$$

with inverse

$$x(t, z) = y(t, z) + \int_0^z L(z, s) y(t, s) ds, \text{ for all } (t, z) \in \Re_+ \times [0,1] \qquad (3.7)$$

maps the solutions of (3.4), (3.5) to the solutions of

$$\frac{\partial y}{\partial t}(t, z) - p \frac{\partial^2 y}{\partial z^2}(t, z) + c\, y(t, z) = 0, \qquad (3.8)$$

$$b_1 y(t,0) + b_2 \frac{\partial y}{\partial z}(t,0) = y(t,1) - u(t) + \int_0^1 k(s) x(t, s) ds = 0, \qquad (3.9)$$

where $k(z) = K(1, z)$ for all $z \in [0,1]$. Selecting $c \geq 0$ sufficiently large, we can guarantee that the closed-loop system (3.4), (3.5), with $u(t) = \int_0^1 k(s) x(t, s) ds$ is exponentially stable in the $L^2$ norm.

Based on the feedback law $u(t) = \int_0^1 k(s) x(t, s) ds$ we are in a position to prove the following result, which shows that the boundary sampled-data feedback law obtained by the emulation of the continuous-time backstepping design works.

**Theorem 3.2:** *Consider the Sturm-Liouville operator $A: D \to L^2([0,1])$ defined by (2.1), (2.2), where $r(z) \equiv 1$, $p(z) \equiv p > 0$, $a_1 = 1$, $a_2 = 0$, $b_1, b_2$ are real constants with $|b_1| + |b_2| > 0$, under Assumption (H). Let $c \geq 0$ be sufficiently large so that the closed-loop system (3.4), (3.5), with $u(t) = \int_0^1 k(s) x(t, s) ds$ is exponentially stable in the $L^2$ norm. Then there exist constants $M, \sigma, T > 0$ such that for every increasing sequence $\{\tau_i \geq 0, i = 0,1,2,...\}$ with $\tau_0 = 0$, $\sup_{i \geq 0}(\tau_{i+1} - \tau_i) \leq T$, $\lim_{i \to +\infty}(\tau_i) = +\infty$ and for every $x_0 \in L^2([0,1])$ the unique function $x: \Re_+ \times [0,1] \to \Re$ for which the mapping $\Re_+ \ni t \to x[t] \in L^2([0,1])$ is continuous, with $x \in C^1(I \times [0,1])$ satisfying $x[t] \in C^2([0,1])$ for all $t > 0$, $x(0, z) = x_0(z)$ for all $z \in [0,1]$, and*

$$\frac{\partial x}{\partial t}(t, z) = p \frac{\partial^2 x}{\partial z^2}(t, z) - q(z) x(t, z), \text{ for all } (t, z) \in I \times (0,1) \qquad (3.10)$$

$$b_1 x(t,0) + b_2 \frac{\partial x}{\partial z}(t,0) = x(t,1) - u(t) = 0, \text{ for all } t \in I \qquad (3.11)$$

$$u(t) = \int_0^1 k(s) x(\tau_i, s) ds, \text{ for } t \in [\tau_i, \tau_{i+1}) \text{ and for all } i \in Z^+ \qquad (3.12)$$

*where $I = \Re_+ \setminus \{\tau_i \geq 0, i = 0,1,2,...\}$, satisfies the following estimate*

$$\|x[t]\|_2 \leq M \exp(-\sigma t) \|x_0\|_2, \text{ for all } t \geq 0. \qquad (3.13)$$

### 3.III. Discussion of the Main Results

Both main results (Theorem 3.1 and Theorem 3.2) provide exponential stability estimates in the $L^2$ norm. Unfortunately, it is not known if exponential stability estimates in different norms hold.

The aforementioned issue is closely related with the consideration of only the case of Dirichlet actuation in the backstepping design. In the cases of Neumann or Robin actuation, the continuous-time backstepping feedback design involves a trace term, i.e., it is of the form $u(t) = k_0 x(t,1) + \int_0^1 k(z) x(t,z) dz$, where $k_0$ is a non-zero real constant and $k \in C^0([0,1])$ is a continuous function. Since the proof of Theorem 3.2 is based on the estimation of the difference between the continuous-time feedback and the applied control action, it follows that in the cases of Neumann or Robin actuation we would have to obtain pointwise estimates for the solution. This is exactly what is needed for the derivation of exponential stability estimates in the $L^\infty$ norm: unfortunately, it is not known if this is possible.

The proof of Theorem 3.1 is very different from the proof of Theorem 3.2. In the proof of Theorem 3.1 there is no need to obtain estimates of the difference between the continuous-time feedback and the applied control action. There is an important reason that explains this difference: the reduced model design and the backstepping design are very different (although they look similar). The backstepping design affects the whole spectrum of the Sturm-Liouville operator $A: D \to C^0([0,1]; \Re)$, while not assigning any of the eigenvalues individually, whereas the reduced model design affects only the first $m$ eigenvalues of $A$, where $m$ is the dimension of the linear finite-dimensional system (3.1). However, both proofs exploit the eigenfunction expansion procedure, which was recently used in [13,15] for the derivation of ISS estimates with respect to boundary disturbances.

It should be noticed that, if the dimension of the linear finite-dimensional system (3.1) $m$ is very large, then the design of a vector $k = (k_1,...,k_m)' \in \Re^m$ for which the matrix $diag(-\lambda_1,...,-\lambda_m) + gk'$ is Hurwitz, can be very difficult. Indeed, if $m$ is very large then the design of $k = (k_1,...,k_m)' \in \Re^m$ can become very sensitive to errors in the constants $g_n$ and the eigenvalues $\lambda_n$ ($n = 1,...,m$); particularly when $\lambda_n < 0$ for $n = 1,...,m$ (strongly unstable open-loop system). This is an important disadvantage of the reduced model design relative to the backstepping design.

Finally, we notice that both main results guarantee robustness with respect to perturbations of the sampling schedule: the exponential stability estimates (3.3) and (3.13) hold for every increasing sequence $\{\tau_i \geq 0, i = 0,1,2,...\}$ with $\tau_0 = 0$, $\sup_{i \geq 0}(\tau_{i+1} - \tau_i) \leq T$, $\lim_{i \to +\infty}(\tau_i) = +\infty$.

## 4. Proofs of Main Results

In order to prove Theorem 3.1, we first need an auxiliary result that guarantees exponential sampled-data stabilization of the linear finite-dimensional system (3.1).

**Lemma 4.1:** *Let* $\lambda_1 < \lambda_2 < ... < \lambda_n < ...$ *with* $\lim_{n \to \infty}(\lambda_n) = +\infty$ *be the eigenvalues of the Sturm-Liouville operator* $A: D \to L_r^2([0,1])$, *defined by (2.1), (2.2). Consider the linear system (3.1), where* $g_n \neq 0$ *($n = 1,...,m$) are real constants,* $g = (g_1,...,g_m)' \in \Re^m$ *and* $x(t) = (x_1(t),...,x_m(t))' \in \Re^m$. *Then the above system is controllable for every integer* $m > 0$. *Moreover, for every vector* $k \in \Re^m$ *for which the matrix* $diag(-\lambda_1,...,-\lambda_m) + rk'$ *is Hurwitz, there exist constants* $M, \sigma, T > 0$ *such that for every increasing sequence* $\{\tau_i \geq 0, i = 0,1,2,...\}$ *with* $\tau_0 = 0$, $\sup_{i \geq 0}(\tau_{i+1} - \tau_i) \leq T$, $\lim_{i \to +\infty}(\tau_i) = +\infty$ *and for every* $x_0 \in \Re^m$ *the solution of the initial value problem (3.1) with*

$$u(t) = k'x(\tau_i), \text{ for } t \in [\tau_i, \tau_{i+1}) \text{ and for all } i \in Z^+ \quad (4.1)$$

*and initial condition* $x(0) = x_0$ *satisfies the following estimate*

$$|x(t)| + |u(t)| \leq M \exp(-\sigma t)|x_0|, \text{ for all } t \geq 0. \tag{4.2}$$

**Proof:** We consider the change of coordinates

$$\xi_n = g_n^{-1} x_n, \quad n = 1,\ldots,m \tag{4.3}$$

We get from (3.1):

$$\dot{\xi} = diag(-\lambda_1,\ldots,-\lambda_m)\xi + bu \tag{4.4}$$

where $b = (1,\ldots,1)' \in \Re^m$. Theorem 3 on page 89 in [32] (Kalman controllability condition) implies that system (4.4) (and consequently its equivalent system (3.1)) is controllable if and only if the rank of the matrix

$$\begin{bmatrix} 1 & -\lambda_1 & \ldots & (-\lambda_1)^{m-1} \\ 1 & -\lambda_2 & \ldots & (-\lambda_2)^{m-1} \\ \vdots & \vdots & & \vdots \\ 1 & -\lambda_m & \ldots & (-\lambda_m)^{m-1} \end{bmatrix}$$

is $n$. However, the above matrix is a Vandermonde matrix and since $\lambda_1 < \lambda_2 < \ldots < \lambda_m$ its determinant is not zero. Thus system (3.1) is controllable.

Let $k \in \Re^m$ be a vector for which the matrix $W := diag(-\lambda_1,\ldots,-\lambda_m) + gk'$ is Hurwitz. Let $T > 0$ be a constant (to be selected). It follows that there exist constants $G, \sigma, \varepsilon > 0$ so that $|\exp(Wt)| \leq G\exp(-(\sigma+\varepsilon)t)$ for all $t \geq 0$. Let $\{\tau_i \geq 0, i = 0,1,2,\ldots\}$ be an increasing sequence with $\tau_0 = 0$, $\sup_{i \geq 0}(\tau_{i+1} - \tau_i) \leq T$, $\lim_{i \to +\infty}(\tau_i) = +\infty$ and let $x_0 \in \Re^m$ be given. Consider the solution $x(t) \in \Re^m$ of the initial value problem (3.1) with (4.1) and initial condition $x(0) = x_0$. It holds that

$$x(t) = \exp(Wt)x(0) + \int_0^t \exp(W(t-s))g(u(s) - k'x(s))ds, \text{ for all } t \geq 0 \tag{4.5}$$

from which we obtain

$$\sup_{0 \leq s \leq t}(|x(s)|\exp(\sigma s)) \leq G|x_0| + G\varepsilon^{-1}|g| \sup_{0 \leq s \leq t}(|u(s) - k'x(s)|\exp(\sigma s)), \text{ for all } t \geq 0 \tag{4.6}$$

For $i \in Z^+$ and $t \in [\tau_i, \tau_{i+1})$ we get (3.1), (4.1) for $n = 1,\ldots,m$:

$$x_n(t) - x_n(\tau_i) = p_n(t - \tau_i)(g_n k'x(\tau_i) - \lambda_n x_n(\tau_i)) = p_n(t - \tau_i)(g_n k - \lambda_n e_n)' x(\tau_i) \tag{4.7}$$

where $e_1 = (1,0,\ldots,0)',\ldots, e_m = (0,0,\ldots,1)'$, $p_n(s) := \dfrac{1 - \exp(-\lambda_n s)}{\lambda_n}$ when $\lambda_n \neq 0$, $s \geq 0$ and $p_n(s) := s$ when $\lambda_n = 0$, $s \geq 0$. Since $p_n(s) \leq p_1(s)$ for all $n = 1,\ldots,m$ and $s \geq 0$, we obtain from (4.7) for all $i \in Z^+$ and $t \in [\tau_i, \tau_{i+1})$:

$$|x(t) - x(\tau_i)| \leq p_1(t - \tau_i)\Gamma|x(\tau_i)| \tag{4.8}$$

where $\Gamma := \left(\sum_{n=1}^m |g_n k - \lambda_n e_n|^2\right)^{1/2}$. Since $p_1(s) \leq p_1(t)$ for all $t \geq s \geq 0$ and since $\sup_{i \geq 0}(\tau_{i+1} - \tau_i) \leq T$, we obtain from (4.1), (4.7) for all $t \geq 0$:

$$|u(t) - k'x(t)| \exp(\sigma t) \leq |k| p_1(T) \exp(\sigma T) \Gamma \sup_{0 \leq s \leq t} (|x(s)| \exp(\sigma s)) \quad (4.9)$$

Combining (4.6) and (4.9), we obtain:

$$\sup_{0 \leq s \leq t} (|x(s)| \exp(\sigma s)) \leq G|x_0| + G\varepsilon^{-1}|g|p_1(T)\exp(\sigma T)|k|\Gamma \sup_{0 \leq s \leq t} (|x(s)| \exp(\sigma s)), \text{ for all } t \geq 0 \quad (4.10)$$

Selecting $T > 0$ sufficiently small so that $G\varepsilon^{-1}|g|p_1(T)\exp(\sigma T)|k|\Gamma < 1$ (always possible since $p_1(s)$ is a continuous function with $p_1(0) = 0$), we conclude that

$$\sup_{0 \leq s \leq t} (|x(s)| \exp(\sigma s)) \leq (1 - G\varepsilon^{-1}|g|p_1(T)\exp(\sigma T)|k|\Gamma) G|x_0|, \text{ for all } t \geq 0 \quad (4.11)$$

Estimate (4.2) with $M = (1 + \exp(\sigma T))(1 - G\varepsilon^{-1}|g|p_1(T)\exp(\sigma T)|k|\Gamma)G$ is a direct consequence of inequality (4.11) and the inequality

$$\sup_{0 \leq s \leq t} (|u(s)| \exp(\sigma s)) \leq \exp(\sigma T) \sup_{0 \leq s \leq t} (|x(s)| \exp(\sigma s)), \text{ for all } t \geq 0$$

which follows from (4.1) and the fact that $\sup_{i \geq 0}(\tau_{i+1} - \tau_i) \leq T$. The proof is complete. ◁

We are now ready to provide the proof of Theorem 3.1.

**Proof of Theorem 3.1:** The fact that $g_n := \frac{p(1)}{a_1^2 + a_2^2}\left(a_2 \phi_n(1) - a_1 \frac{d\phi_n}{dz}(1)\right) \neq 0$ ($n = 1,\ldots,m$) are not zero is a direct consequence of the boundary condition $a_1 \phi_n(1) + a_2 \frac{d\phi_n}{dz}(1) = 0$ (otherwise we would have $\phi_n(1) = \frac{d\phi_n}{dz}(1) = 0$ and consequently, $\phi_n(z) \equiv 0$). Let $M, \sigma, T > 0$ be constants such that for every increasing sequence $\{\tau_i \geq 0, i = 0,1,2,\ldots\}$ with $\tau_0 = 0$, $\sup_{i \geq 0}(\tau_{i+1} - \tau_i) \leq T$, $\lim_{i \to +\infty}(\tau_i) = +\infty$ and for every $x_0 \in \Re^m$ the solution of the initial value problem (3.1) with (4.1) and initial condition $x(0) = x_0$ satisfies estimate (4.2).

Let $\{\tau_i \geq 0, i = 0,1,2,\ldots\}$ be an increasing sequence with $\tau_0 = 0$, $\sup_{i \geq 0}(\tau_{i+1} - \tau_i) \leq T$, $\lim_{i \to +\infty}(\tau_i) = +\infty$ and let $x_0 \in L_r^2([0,1])$ be an arbitrary function.

Existence/uniqueness of the function $x: \Re_+ \times [0,1] \to \Re$ for which the mapping $\Re_+ \ni t \to x[t] \in L_r^2([0,1])$ is continuous, with $x \in C^1(I \times [0,1])$ satisfying $x[t] \in C^2([0,1])$ for all $t > 0$, $x(0,z) = x_0(z)$ for all $z \in [0,1]$, (2.9), (2.10) and (3.2) is guaranteed by Corollary 2.2.

Since the eigenfunctions $\{\phi_n\}_{n=1}^{\infty}$ of the Sturm-Liouville operator $A: D \to L_r^2([0,1])$ defined by (2.1), (2.2) under Assumption (H) form an orthonormal basis of $L_r^2([0,1])$, it follows that Parseval's identity holds, i.e.,

$$\|x[t]\|_r^2 = \sum_{n=1}^{\infty} x_n^2(t), \text{ for all } t \geq 0 \quad (4.12)$$

where

$$x_n(t) := \int_0^1 r(z) x(t,z) \phi_n(z) dz, \text{ for } n = 1,2,\ldots \quad (4.13)$$

Notice that the mappings $\mathfrak{R}_+ \ni t \to x_n(t)$ for $n=1,2,...$ are continuous (by virtue of continuity of the mapping $\mathfrak{R}_+ \ni t \to x[t] \in L_r^2([0,1])$) and the mappings $I \ni t \to x_n(t)$ for $n=1,2,...$ are $C^1$. By virtue of (2.9), it follows from repeated integration by parts, that the following equalities hold for all $t \in I$:

$$\dot{x}_n(t) = \int_0^1 r(z) \frac{\partial x}{\partial t}(t,z) \phi_n(z) dz = \int_0^1 \frac{\partial}{\partial z}\left( p(z) \frac{\partial x}{\partial z}(t,z) \right) \phi_n(z) dz - \int_0^1 q(z) x(t,z) \phi_n(z) dz$$

$$= p(1) \frac{\partial x}{\partial z}(t,1) \phi_n(1) - p(0) \frac{\partial x}{\partial z}(t,0) \phi_n(0) - \int_0^1 p(z) \frac{\partial x}{\partial z}(t,z) \frac{d\phi_n}{dz}(z) dz - \int_0^1 q(z) x(t,z) \phi_n(z) dz =$$

$$= p(1)\left( \frac{\partial x}{\partial z}(t,1) \phi_n(1) - x(t,1) \frac{d\phi_n}{dz}(1) \right) + p(0)\left( \frac{d\phi_n}{dz}(0) x(t,0) - \frac{\partial x}{\partial z}(t,0) \phi_n(0) \right)$$

$$+ \int_0^1 x(t,z) \left[ \frac{d}{dz}\left( p(z) \frac{d\phi_n}{dz}(z) \right) - q(z) \phi_n(z) \right] dz$$

Thus we get for all $t \in I$:

$$\dot{x}_n(t) = p(1)\left( \frac{\partial x}{\partial z}(t,1) \phi_n(1) - x(t,1) \frac{d\phi_n}{dz}(1) \right) + p(0)\left( \frac{d\phi_n}{dz}(0) x(t,0) - \frac{\partial x}{\partial z}(t,0) \phi_n(0) \right)$$

$$- \int_0^1 r(z) x(t,z) (A\phi_n)(z) dz \tag{4.14}$$

It follows from (4.14), the fact that $(A\phi_n)(z) = \lambda_n \phi_n(z)$ and definition (4.13) that the following equation holds for all $t \in I$:

$$\dot{x}_n(t) + \lambda_n x_n(t) = p(1)\left( \frac{\partial x}{\partial z}(t,1) \phi_n(1) - x(t,1) \frac{d\phi_n}{dz}(1) \right)$$

$$+ p(0)\left( \frac{d\phi_n}{dz}(0) x(t,0) - \frac{\partial x}{\partial z}(t,0) \phi_n(0) \right) \tag{4.15}$$

Next, we show that for all $t \in I$:

$$x(t,0) \frac{d\phi_n}{dz}(0) - \frac{\partial x}{\partial z}(t,0) \phi_n(0) = 0 \tag{4.16}$$

Indeed, the equation $b_1 x(t,0) + b_2 \frac{\partial x}{\partial z}(t,0) = 0$ (see (2.10)) gives:

$$b_1^2 x(t,0) \frac{d\phi_n}{dz}(0) + b_1 b_2 \frac{d\phi_n}{dz}(0) \frac{\partial x}{\partial z}(t,0) = 0$$

$$b_1 b_2 \phi_n(0) x(t,0) + b_2^2 \phi_n(0) \frac{\partial x}{\partial z}(t,0) = 0$$

from which we obtain:

$$0 = (b_1^2 + b_2^2)\left( x(t,0) \frac{d\phi_n}{dz}(0) - \phi_n(0) \frac{\partial x}{\partial z}(t,0) \right)$$

$$+ \left( b_1 \frac{\partial x}{\partial z}(t,0) - b_2 x(t,0) \right)\left( b_2 \frac{d\phi_n}{dz}(0) + b_1 \phi_n(0) \right)$$

Equation (4.16) follows directly from the above equation and the fact that $b_1 \phi_n(0) + b_2 \frac{d\phi_n}{dz}(0) = 0$.

Moreover, the equation $a_1 x(t,1) + a_2 \frac{\partial x}{\partial z}(t,1) = u(t)$ (see (2.10)) gives:

$$a_1^2 x(t,1) \frac{d\phi_n}{dz}(1) + a_1 a_2 \frac{d\phi_n}{dz}(1) \frac{\partial x}{\partial z}(t,1) = u(t) a_1 \frac{d\phi_n}{dz}(1)$$

$$a_1 a_2 \phi_n(1) x(t,1) + a_2^2 \phi_n(1) \frac{\partial x}{\partial z}(t,1) = u(t) a_2 \phi_n(1)$$

from which we obtain:

$$u(t)\left(a_1 \frac{d\phi_n}{dz}(1) - a_2\phi_n(1)\right) = (a_1^2 + a_2^2)\left(x(t,1)\frac{d\phi_n}{dz}(1) - \phi_n(1)\frac{\partial x}{\partial z}(t,1)\right)$$
$$+ \left(a_1 \frac{\partial x}{\partial z}(t,1) - a_2 x(t,1)\right)\left(a_2 \frac{d\phi_n}{dz}(1) + a_1\phi_n(1)\right)$$

The fact that $a_1\phi_n(1) + a_2 \frac{d\phi_n}{dz}(1) = 0$ in conjunction with the above equation implies that:

$$\phi_n(1)\frac{\partial x}{\partial z}(t,1) - x(t,1)\frac{d\phi_n}{dz}(1) = \frac{u(t)}{a_1^2 + a_2^2}\left(a_2\phi_n(1) - a_1 \frac{d\phi_n}{dz}(1)\right) \quad (4.17)$$

Using (4.15), (4.16) and (4.17), we obtain for all $t \in I$:

$$\dot{x}_n(t) + \lambda_n x_n(t) = \frac{p(1)u(t)}{a_1^2 + a_2^2}\left(a_2\phi_n(1) - a_1 \frac{d\phi_n}{dz}(1)\right) \quad (4.18)$$

It follows from (4.2) that the following inequality holds:

$$|u(t)| + \left(\sum_{n=1}^{m} x_n^2(t)\right)^{1/2} \leq M \exp(-\sigma t)\left(\sum_{n=1}^{m} x_n^2(0)\right)^{1/2}, \text{ for all } t \geq 0. \quad (4.19)$$

Without loss of generality we may assume that $\sigma < \lambda_{m+1}$. Integrating the differential equations (4.18), we obtain for all $t \geq 0$ and $n = m+1, m+2,...$:

$$x_n(t) = \exp(-\lambda_n t)x_n(0) + \frac{p(1)}{a_1^2 + a_2^2}\left(a_2\phi_n(1) - a_1 \frac{d\phi_n}{dz}(1)\right)\int_0^t \exp(-\lambda_n(t-s))u(s)ds \quad (4.20)$$

Equations (4.20) in conjunction with the inequality $\sigma < \lambda_{m+1}$ imply the following estimates for all $t \geq 0$ and $n = m+1, m+2,...$:

$$|x_n(t)| \leq \exp(-\sigma t)|x_n(0)| + \frac{p(1)}{a_1^2 + a_2^2}\left|a_2\phi_n(1) - a_1 \frac{d\phi_n}{dz}(1)\right|\frac{1}{\lambda_n - \sigma}\max_{0 \leq s \leq t}(|u(s)|\exp(-\sigma(t-s))) \quad (4.21)$$

Let $w > -\lambda_1$ be a constant. Notice that the Sturm-Liouville operator $(A + wId)$, where $Id$ is the identity operator, has positive eigenvalues ($\lambda_n + w$ for $n = 1,2,...$) and the same eigenfunctions with $A$. Changing $z$ by $1-z$ and using Lemma 2.1 in [13], it follows that the boundary value problem

$$\frac{d}{dz}\left(p(z)\frac{d\bar{x}}{dz}(z)\right) - q(z)\bar{x}(z) - w\bar{x}(z) = 0, \text{ for all } z \in [0,1], \quad (4.22)$$

with

$$b_1\bar{x}(0) + b_2 \frac{d\bar{x}}{dz}(0) = 0, \quad a_1\bar{x}(1) + a_2 \frac{d\bar{x}}{dz}(1) = \sqrt{a_1^2 + a_2^2} \quad (4.23)$$

has a unique solution $\bar{x} \in C^2([0,1])$, which satisfies

$$p^2(1)\sum_{n=1}^{\infty}(\lambda_n + w)^{-2}\left|\frac{a_1}{\sqrt{a_1^2 + a_2^2}}\frac{d\phi_n}{dz}(1) - \frac{a_2}{\sqrt{a_1^2 + a_2^2}}\phi_n(1)\right|^2 = \int_0^1 r(z)\bar{x}^2(z)dz \quad (4.24)$$

Consequently, we obtain from (4.24):

$$K := p^2(1)\sum_{n=m+1}^{\infty}(\lambda_n - \sigma)^{-2}\left|\frac{a_1}{\sqrt{a_1^2 + a_2^2}}\frac{d\phi_n}{dz}(1) - \frac{a_2}{\sqrt{a_1^2 + a_2^2}}\phi_n(1)\right|^2 \leq \left(\frac{\lambda_{m+1} + w}{\lambda_{m+1} - \sigma}\right)^2 \int_0^1 r(z)\bar{x}^2(z)dz < +\infty \quad (4.25)$$

It follows from (4.21) and (4.25) that the following estimate holds for all $t \geq 0$:

$$\sum_{n=m+1}^{\infty}|x_n(t)|^2 \leq 2\exp(-2\sigma t)\sum_{n=m+1}^{\infty}|x_n(0)|^2 + \frac{2K}{a_1^2+a_2^2}\max_{0\leq s\leq t}\left(|u(s)|^2\exp(-2\sigma(t-s))\right) \quad (4.26)$$

Estimate (3.3) is a direct consequence of (4.12), (4.19) and (4.26). The proof is complete. ◁

Next, we provide the proof of Theorem 3.2.

**Proof of Theorem 3.2:** Let a constant $T>0$ be a constant (to be selected). Let $\{\tau_i \geq 0, i=0,1,2,...\}$ be an increasing sequence with $\tau_0 = 0$, $\sup_{i\geq 0}(\tau_{i+1}-\tau_i)\leq T$, $\lim_{i\to+\infty}(\tau_i)=+\infty$ and let $x_0 \in L^2([0,1])$ be given. Existence/uniqueness of the function $x:\Re_+\times[0,1]\to\Re$ for which the mapping $\Re_+ \ni t \to x[t]\in L^2([0,1])$ is continuous, with $x\in C^1(I\times[0,1])$ satisfying $x[t]\in C^2([0,1])$ for all $t>0$, $x(0,z)=x_0(z)$ for all $z\in[0,1]$, (3.10), (3.11) and (3.12) is guaranteed by Corollary 2.2. It follows from (3.6) that the function $y:\Re_+\times[0,1]\to\Re$ defined by (3.6) is of class $y\in C^1(I\times[0,1])$ satisfying $y[t]\in C^2([0,1])$ for all $t>0$, $y(0,z)=x_0(z)-\int_0^z K(z,s)x_0(s)ds$ for all $z\in[0,1]$, (3.8) for all $(t,z)\in I\times(0,1)$ and (3.9) for all $t\in I$.

Let $0<\mu_1<\mu_2<...<\mu_n<...$ with $\lim_{n\to\infty}(\mu_n)=+\infty$ be the eigenvalues of the Sturm-Liouville operator $B:D\to L^2([0,1])$ defined by (2.1), (2.2) with $r(z)\equiv 1$, $p(z)\equiv p$, $q(z)\equiv c$, $a_1=1$, $a_2=0$. Let $\{\psi_n\}_{n=1}^{\infty}$ be the eigenfunctions of the operator $B:D\to L^2([0,1])$. It should be noticed that, as remarked in [13] the operator $B:D\to L^2([0,1])$ satisfies Assumption (H). Since the eigenfunctions of the Sturm-Liouville operator $B:D\to L^2([0,1])$ form an orthonormal basis of $L^2([0,1])$, it follows that Parseval's identity holds, i.e.,

$$\|y[t]\|_2^2 = \sum_{n=1}^{\infty} y_n^2(t), \text{ for all } t\geq 0 \quad (4.27)$$

where

$$y_n(t) := \int_0^1 y(t,z)\psi_n(z)dz, \text{ for } n=1,2,... \quad (4.28)$$

Notice that the mappings $\Re_+ \ni t \to y_n(t)$ for $n=1,2,...$ are continuous (by virtue of continuity of the mapping $\Re_+ \ni t \to x[t]\in L^2([0,1])$ and by virtue of continuity of the transformation (3.6)) and the mappings $I \ni t \to y_n(t)$ for $n=1,2,...$ are $C^1$. Following the same procedure as in the proof of Theorem 3.1, we are in a position to guarantee that the following equalities hold for all $t\in I$:

$$\dot{y}_n(t) + \mu_n y_n(t) = -p\frac{d\psi_n}{dz}(1)v(t) \quad (4.29)$$

where

$$v(t) := u(t) - \int_0^1 k(s)x(t,s)ds, \text{ for all } t\geq 0 \quad (4.30)$$

Following the same procedure as in the proof of Theorem 2.2 in [13] and using the continuity of the mappings $\Re_+ \ni t \to y_n(t)$ for $n=1,2,...$ as well as (4.27) and (4.29), we are in a position to guarantee the existence of constants $G,\sigma,\gamma>0$ such that the following estimate holds for all $t\geq 0$:

$$\|y[t]\|_2 \leq G\exp(-\sigma t)\|y[0]\|_2 + \gamma\sup_{0\leq s\leq t}\left(|v(s)|\exp(-\sigma(t-s))\right) \quad (4.31)$$

A detailed derivation of (4.31) is given in the Appendix.

Next, define the following upper bounds for the norms of the transformations given by (3.6), (3.7):

$$\tilde{K} := 1 + \left( \int_0^1 \left( \int_0^z |K(z,s)|^2 ds \right) dz \right)^{1/2} \quad , \quad \tilde{L} := 1 + \left( \int_0^1 \left( \int_0^z |L(z,s)|^2 ds \right) dz \right)^{1/2} \quad (4.32)$$

Moreover, let $N \geq 1$ be an integer (to be selected) and define:

$$g(z) := \sum_{n=1}^{N} k_n \phi_n(z), \text{ for all } z \in [0,1] \quad (4.33)$$

where

$$k_n := \int_0^1 k(s) \phi_n(s) ds \quad (4.34)$$

We select $N \geq 1$ to be sufficiently large so that

$$2\gamma \tilde{L} \|k - g\|_2 < 1 \quad (4.35)$$

Define:

$$w(t) = \int_0^1 g(s)(x(\tau_i, s) - x(t, s)) ds, \text{ for } t \in [\tau_i, \tau_{i+1}) \text{ and for all } i \in Z^+ \quad (4.36)$$

Notice that (3.12) and definitions (4.30), (4.36) imply that

$$v(t) = w(t) + \int_0^1 (k(s) - g(s))(x(\tau_i, s) - x(t, s)) ds, \text{ for } t \in [\tau_i, \tau_{i+1}) \text{ and for all } i \in Z^+ \quad (4.37)$$

Using (3.10) we obtain for $t \in (\tau_i, \tau_{i+1})$ and for all $i \in Z^+$:

$$\dot{w}(t) = p\left( \frac{dg}{dz}(1) x(t,1) - g(1) \frac{\partial x}{\partial z}(t,1) \right) + p\left( g(0) \frac{\partial x}{\partial z}(t,0) - \frac{dg}{dz}(0) x(t,0) \right) + \int_0^1 (Ag)(s) x(t,s) ds \quad (4.38)$$

Using (4.16), (3.11), (3.12), definition (4.33) and the fact that $\phi_n(1) = 0$ for $n = 1,2,...$, we obtain from (4.38) for $t \in (\tau_i, \tau_{i+1})$ and for all $i \in Z^+$:

$$\dot{w}(t) = p \int_0^1 k(s) x(\tau_i, s) ds \sum_{n=1}^{N} k_n \frac{d\phi_n}{dz}(1) + \sum_{n=1}^{N} k_n \lambda_n \int_0^1 \phi_n(s) x(t,s) ds \quad (4.39)$$

It follows from (4.39), the Cauchy-Schwarz inequality and the fact that $\|\phi_n\| = 1$ for $n = 1,2,...$ that the following inequality holds $t \in (\tau_i, \tau_{i+1})$ and for all $i \in Z^+$:

$$|\dot{w}(t)| \leq p \|k\|_2 \|x[\tau_i]\|_2 \sum_{n=1}^{N} \left| k_n \frac{d\phi_n}{dz}(1) \right| + \|x[t]\|_2 \sum_{n=1}^{N} |k_n \lambda_n| \quad (4.40)$$

Continuity of the mapping $\Re_+ \ni t \to x[t] \in L^2([0,1])$ in conjunction with (4.36) implies that the mapping $\Re_+ \ni t \to w(t) \in \Re$ is right continuous with $w(\tau_i) = 0$ for all $i \in Z^+$. Consequently, we get from (4.40) in conjunction with the fact that $\sup_{i \geq 0}(\tau_{i+1} - \tau_i) \leq T$:

$$|w(t)| \leq Tp \|k\|_2 \|x[\tau_i]\|_2 \sum_{n=1}^{N} \left| k_n \frac{d\phi_n}{dz}(1) \right| + T \sup_{\tau_i \leq s \leq t} (\|x[t]\|_2) \sum_{n=1}^{N} |k_n \lambda_n|, \text{ for } t \in [\tau_i, \tau_{i+1}) \text{ and for all } i \in Z^+ \quad (4.41)$$

Using (4.41) and the fact that $\sup_{i \geq 0}(\tau_{i+1} - \tau_i) \leq T$, we obtain:

$$|w(t)|\exp(\sigma t) \leq T\exp(\sigma T)\left(p\|k\|_2 \sum_{n=1}^{N}\left|k_n \frac{d\phi_n}{dz}(1)\right| + \sum_{n=1}^{N}|k_n \lambda_n|\right)\sup_{0\leq s\leq t}\left(\|x[s]\|_2 \exp(\sigma s)\right), \text{ for } t\geq 0 \quad (4.42)$$

Using (4.37) and (4.42), we obtain for all $t \in [\tau_i, \tau_{i+1})$ and for all $i \in Z^+$:

$$|v(t)|\exp(\sigma t) \leq T\exp(\sigma T)\left(p\|k\|_2 \sum_{n=1}^{N}\left|k_n \frac{d\phi_n}{dz}(1)\right| + \sum_{n=1}^{N}|k_n \lambda_n|\right)\sup_{0\leq s\leq t}\left(\|x[s]\|_2 \exp(\sigma s)\right)$$
$$+\|k-g\|_2 \|x[\tau_i]\|_2 \exp(\sigma t) + \|k-g\|_2 \|x[t]\|_2 \exp(\sigma t) \quad (4.43)$$

It follows from (4.43) and the fact that $\sup_{i\geq 0}(\tau_{i+1}-\tau_i)\leq T$ that the following estimate holds for all $t\geq 0$:

$$|v(t)|\exp(\sigma t) \leq T\exp(\sigma T)\left(p\|k\|_2 \sum_{n=1}^{N}\left|k_n \frac{d\phi_n}{dz}(1)\right| + \sum_{n=1}^{N}|k_n \lambda_n|\right)\sup_{0\leq s\leq t}\left(\|x[s]\|_2 \exp(\sigma s)\right)$$
$$+\|k-g\|_2 (\exp(\sigma T)+1)\sup_{0\leq s\leq t}\left(\|x[s]\|_2 \exp(\sigma s)\right) \quad (4.44)$$

Using the fact that $\tilde{L}>0$ defined by (4.32) is an upper bound for the norm of the transformation (3.7), we obtain for all $t\geq 0$:

$$|v(t)|\exp(\sigma t) \leq \tilde{L}T\exp(\sigma T)\left(p\|k\|_2 \sum_{n=1}^{N}\left|k_n \frac{d\phi_n}{dz}(1)\right| + \sum_{n=1}^{N}|k_n \lambda_n|\right)\sup_{0\leq s\leq t}\left(\|y[s]\|_2 \exp(\sigma s)\right)$$
$$+\tilde{L}\|k-g\|_2 (\exp(\sigma T)+1)\sup_{0\leq s\leq t}\left(\|y[s]\|_2 \exp(\sigma s)\right) \quad (4.45)$$

Combining (4.31) with (4.45), we obtain for all $t\geq 0$:

$$\sup_{0\leq s\leq t}\left(\|y[s]\|_2 \exp(\sigma s)\right) \leq \gamma\tilde{L}T\exp(\sigma T)\left(p\|k\|_2 \sum_{n=1}^{N}\left|k_n \frac{d\phi_n}{dz}(1)\right| + \sum_{n=1}^{N}|k_n \lambda_n|\right)\sup_{0\leq s\leq t}\left(\|y[s]\|_2 \exp(\sigma s)\right)$$
$$+\gamma\tilde{L}\|k-g\|_2 (\exp(\sigma T)+1)\sup_{0\leq s\leq t}\left(\|y[s]\|_2 \exp(\sigma s)\right) + G\|y[0]\|_2 \quad (4.46)$$

We select $T>0$ sufficiently small so that:

$$\gamma\tilde{L}T\exp(\sigma T)\left(p\|k\|_2 \sum_{n=1}^{N}\left|k_n \frac{d\phi_n}{dz}(1)\right| + \sum_{n=1}^{N}|k_n \lambda_n|\right) + \gamma\tilde{L}\|k-g\|_2 (\exp(\sigma T)+1) < 1 \quad (4.47)$$

The existence of a constant $T>0$ sufficiently small so that (4.47) holds is a consequence of the continuity of the expression of the left hand side of inequality (4.47) with respect to $T$ and inequality (4.35). Inequality (4.46) in conjunction with (4.47) and the fact that $\tilde{K}>0$ defined by (4.32) is an upper bound for the norm of the transformation (3.6), gives for all $t\geq 0$:

$$\sup_{0\leq s\leq t}\left(\|y[s]\|_2 \exp(\sigma s)\right)$$
$$\leq \left(1-\gamma\tilde{L}T\exp(\sigma T)\left(p\|k\|_2 \sum_{n=1}^{N}\left|k_n \frac{d\phi_n}{dz}(1)\right| + \sum_{n=1}^{N}|k_n \lambda_n|\right) - \gamma\tilde{L}\|k-g\|_2 (\exp(\sigma T)+1)\right)^{-1} G\tilde{K}\|x[0]\|_2 \quad (4.48)$$

Inequality (3.13) with

$$M := \left(1-\gamma\tilde{L}T\exp(\sigma T)\left(p\|k\|_2 \sum_{n=1}^{N}\left|k_n \frac{d\phi_n}{dz}(1)\right| + \sum_{n=1}^{N}|k_n \lambda_n|\right) - \gamma\tilde{L}\|k-g\|_2 (\exp(\sigma T)+1)\right)^{-1} G\tilde{L}\tilde{K}$$

is a direct consequence of (4.48) and the fact that $\tilde{L}>0$ defined by (4.32) is an upper bound for the norm of the transformation (3.7). The proof is complete. ◁

**Discussion of the proof of Theorem 3.2:** There are two things that are important in the proof of Theorem 3.2.

**1)** The stability analysis is performed for the transformed system (3.8), (3.9) and not for the original system. This feature is expected since it holds also for the case of the continuous-time feedback (see [16,33]). However, the stability analysis of the transformed system (3.8), (3.9) becomes more involved here, because of the existence of a perturbation in the boundary condition: the difference between the applied control action and the action determined by the continuous-time nominal feedback. This difference is defined in (4.30) to be the signal $v(t)$.

**2)** It can be shown that the time derivative of the signal $v(t)$ satisfies the differential equation

$$\dot{v}(t) = p\left(\frac{dk}{dz}(1)x(t,1) - k(1)\frac{\partial x}{\partial z}(t,1)\right) + p\left(k(0)\frac{\partial x}{\partial z}(t,0) - \frac{dk}{dz}(0)x(t,0)\right) + \int_0^1 (Ak)(s)x(t,s)ds$$

for all $t \in (\tau_i, \tau_{i+1})$ and for all $i \in Z^+$. However, the above differential equation does not allow the derivation of an upper bound of the magnitude of $\dot{v}(t)$, since it contains terms which cannot be estimated (e.g., the term $k(1)\frac{\partial x}{\partial z}(t,1)$). Therefore, in order to estimate the magnitude of the signal $v(t)$, we use the decomposition (4.37). In this way, the magnitude of the signal $w(t)$ can be estimated by estimating the magnitude of its time derivative and the signal $\int_0^1 (k(s) - g(s))(x(\tau_i, s) - x(t,s))ds$ can become arbitrarily small by controlling the magnitude of $\|k - g\|_2$.

## 5. Illustrative Example

We consider the following control system

$$\frac{\partial x}{\partial t}(t,z) = p\frac{\partial^2 x}{\partial z^2}(t,z) + qx(t,z) \tag{5.1}$$

$$x(t,1) = u(t) \tag{5.2}$$

$$x(t,0) = 0, \tag{5.3}$$

where $p > 0, q \in \Re$ are constants, $x[t]$ is the state and $u(t)$ is the control input. More specifically, we consider the system under boundary sampled-data control with ZOH:

$$u(t) = u_i, \text{ for } t \in [\tau_i, \tau_{i+1}) \text{ and for all } i \in Z^+ \tag{5.4}$$

where $\{\tau_i \geq 0, i = 0,1,2,...\}$ is an increasing sequence (the sequence of sampling times) with $\tau_0 = 0$, $\lim_{i \to +\infty}(\tau_i) = +\infty$ and $\{u_i \in \Re, i = 0,1,2,...\}$ is the sequence of applied inputs.

Notice that the Sturm-Liouville operator $A = -p\frac{d^2}{dz^2} - q$ defined on the set of all functions $f \in H^2([0,1])$ for which

$$f(0) = f(1) = 0 \tag{5.5}$$

satisfies Assumption (H) with $r(z) \equiv 1$, $p(z) \equiv p$, $q(z) \equiv -q$, $a_1 = b_1 = 1$, $a_2 = b_2 = 0$, $\phi_n(z) = \sqrt{2}\sin(n\pi z)$, $\lambda_n = n^2\pi^2 p - q$ for $n = 1,2,...$. Suppose that

$$p\pi^2 \leq q < 4p\pi^2 \tag{5.6}$$

Let $k > q - p\pi^2$ be a constant (arbitrary). Applying Theorem 3.1 with $m = 1$, $g_1 := p\sqrt{2}$, we conclude that there exist constants $G, \mu, T > 0$ such that for every increasing sequence $\{\tau_i \geq 0, i = 0,1,2,...\}$ with $\tau_0 = 0$, $\sup_{i \geq 0}(\tau_{i+1} - \tau_i) \leq T$, $\lim_{i \to +\infty}(\tau_i) = +\infty$ and for every $x_0 \in L^2([0,1])$ the unique function $x : \Re_+ \times [0,1] \to \Re$ for which the mapping $\Re_+ \ni t \to x[t] \in L^2([0,1])$ is continuous, with $x \in C^1(I \times [0,1])$ satisfying $x[t] \in C^2([0,1])$ for all $t > 0$, $x(0,z) = x_0(z)$ for all $z \in [0,1]$, and

$$\frac{\partial x}{\partial t}(t,z) = p\frac{\partial^2 x}{\partial z^2}(t,z) + qx(t,z), \text{ for all } (t,z) \in I \times (0,1) \tag{5.7}$$

$$x(t,0) = x(t,1) - u(t) = 0, \text{ for all } t \in I \tag{5.8}$$

$$u(t) = -kp^{-1}\int_0^1 \sin(\pi z)x(\tau_i, z)dz, \text{ for } t \in [\tau_i, \tau_{i+1}) \text{ and for all } i \in Z^+ \tag{5.9}$$

where $I = \Re_+ \setminus \{\tau_i \geq 0, i = 0,1,2,...\}$, satisfies the following estimate

$$\|x[t]\|_2 \leq G\exp(-\mu t)\|x_0\|_2, \text{ for all } t \geq 0. \tag{5.10}$$

Following the proof of Lemma 4.1 and the proof of Theorem 3.1, we are in a position to determine a (conservative) upper bound for $T > 0$ and the convergence rate $\mu > 0$. The constant $T > 0$ must satisfy the inequality

$$k(k + p\pi^2 - q)T\exp((q + \sigma - p\pi^2)T) + \sigma < k + p\pi^2 - q$$

for certain $\sigma > 0$ and the convergence rate $\mu > 0$ satisfies the (conservative) bounds $\mu < 4p\pi^2 - q$ and $\mu \leq \sigma$. Similarly, we can deal with cases where $q \geq 4p\pi^2$.

On the other hand, the results in Section VIII.A. in [33] and Theorem 3.2 allow us to obtain a different sampled-data feedback law for system (5.1), (5.2), (5.3). Let $I_1(s)$ denote the modified Bessel function of order 1. It follows from Theorem 3.2 that for every $c \geq 0$, there exist constants $G, \mu, T > 0$ such that for every increasing sequence $\{\tau_i \geq 0, i = 0,1,2,...\}$ with $\tau_0 = 0$, $\sup_{i \geq 0}(\tau_{i+1} - \tau_i) \leq T$, $\lim_{i \to +\infty}(\tau_i) = +\infty$ and for every $x_0 \in L^2([0,1])$ the unique function $x : \Re_+ \times [0,1] \to \Re$ for which the mapping $\Re_+ \ni t \to x[t] \in L^2([0,1])$ is continuous, with $x \in C^1(I \times [0,1])$ satisfying $x[t] \in C^2([0,1])$ for all $t > 0$, $x(0,z) = x_0(z)$ for all $z \in [0,1]$, (5.7), (5.8) and

$$u(t) = -(q+c)p^{-1}\int_0^1 s\frac{I_1\left(\sqrt{(q+c)p^{-1}(1-s^2)}\right)}{\sqrt{(q+c)p^{-1}(1-s^2)}}x(\tau_i, s)ds, \text{ for } t \in [\tau_i, \tau_{i+1}) \text{ and for all } i \in Z^+ \tag{5.11}$$

where $I = \Re_+ \setminus \{\tau_i \geq 0, i = 0,1,2,...\}$, satisfies (5.10).

It should be noticed that the feedback law (5.11) is guaranteed to work for all values of the constants $p > 0, q \in \Re$, whereas the control law given by (4.8) is only applicable for $p > 0, q \in \Re$ which satisfy inequality (5.6). ◁

# 6. Concluding Remarks

The paper provides two different results for the application of boundary feedback control with ZOH to 1-D linear parabolic systems on bounded domains. The two different results are developed for two different continuous-time boundary feedback designs: the reduced model design and the backstepping design. It is shown that the continuous-time boundary feedback applied in a sample-and-hold fashion guarantees closed-loop exponential stability, provided that the sampling period is sufficiently small. The obtained results provide stability estimates for weighted 2-norms of the state and robustness with respect to perturbations of the sampling schedule is guaranteed.

Future work may involve the development of boundary feedback designs that are capable to handle simultaneous sampling in space and time. To this purpose, sampled-data observers for linear 1-D parabolic systems must be developed.

# Appendix

**Proof of Theorem 2.1:** Without loss of generality we will assume that $\lambda_1 > 0$. Indeed, if this is not the case, we may perform the same analysis for the function $y(t,z) = \exp(-kt)x(t,z)$ with $k > \lambda_1$ (the function $y(t,z)$ satisfies a PDE similar to (2.4) with the corresponding Sturm-Liouville operator satisfying $\lambda_1 > 0$). Moreover, it follows from (2.3) that

$$\sum_{n=1}^{\infty} \lambda_n^{-1} \max_{0 \le z \le 1}(|\phi_n(z)|) < +\infty \tag{A.1}$$

Define

$$v(t) := u(t), \text{ for } t \in [0,T) \tag{A.2}$$

$$v(t) = \lim_{l \to T^-}(u(l)) + (t-T)\lim_{l \to T^-}(\dot{u}(l)) + \frac{1}{2}(t-T)^2 \lim_{l \to T^-}(\ddot{u}(l)), \text{ for } t \in [T, +\infty) \tag{A.3}$$

and consider the solution of

$$\frac{\partial y}{\partial t}(t,z) - \frac{1}{r(z)}\frac{\partial}{\partial z}\left(p(z)\frac{\partial y}{\partial z}(t,z)\right) + \frac{q(z)}{r(z)}y(t,z) = f(t,z), \text{ for all } (t,z) \in (0,+\infty) \times (0,1) \tag{A.4}$$

$$b_1 y(t,0) + b_2 \frac{\partial y}{\partial z}(t,0) = a_1 y(t,1) + a_2 \frac{\partial y}{\partial z}(t,1) = 0, \text{ for all } t \in (0,+\infty) \tag{A.5}$$

where

$$f(0,z) := \frac{v(0)}{r(z)}\left(\frac{d}{dz}\left((2\sigma_1 z + 3\sigma_2 z^2)p(z)\right) - q(z)(\sigma_1 z^2 + \sigma_2 z^3)\right) - (\sigma_1 z^2 + \sigma_2 z^3)\lim_{l \to 0^+}(\dot{u}(l)),$$

$$\text{for all } z \in [0,1] \tag{A.6}$$

$$f(t,z) := \frac{v(t)}{r(z)}\left(\frac{d}{dz}\left((2\sigma_1 z + 3\sigma_2 z^2)p(z)\right) - q(z)(\sigma_1 z^2 + \sigma_2 z^3)\right) - (\sigma_1 z^2 + \sigma_2 z^3)\dot{v}(t),$$

$$\text{for all } (t,z) \in (0,+\infty) \times [0,1] \tag{A.7}$$

where $\sigma_1 := \frac{3a_1 - a_2}{a_1^2 + a_2^2}$, $\sigma_2 := \frac{a_2 - 2a_1}{a_1^2 + a_2^2}$, with initial condition $y(0,z) = x_0(z) - u(0)(\sigma_1 z^2 + \sigma_2 z^3)$ for all $z \in [0,1]$. Notice that $f \in C^1(\Re_+ \times [0,1])$. We will show next that there is a unique function $y : \Re_+ \times [0,1] \to \Re$ for which the mapping $\Re_+ \ni t \to y[t] \in L_r^2([0,1])$ is continuous, with $y \in C^1((0,+\infty) \times [0,1])$ satisfying $y[t] \in C^2([0,1])$ for all $t > 0$, $y(0,z) = x_0(z) - u(0)(\sigma_1 z^2 + \sigma_2 z^3)$ for all $z \in [0,1]$ and (A.4), (A.5). In this way, we define

$$x(t,z) := y(t,z) + v(t)(\sigma_1 z^2 + \sigma_2 z^3), \text{ for all } (t,z) \in (0,+\infty) \times [0,1] \tag{A.8}$$

$$x(0,z) := x_0(z), \text{ for all } z \in [0,1] \tag{A.9}$$

Notice that (A.2) in conjunction with the fact that $y(0,z) = x_0(z) - u(0)(\sigma_1 z^2 + \sigma_2 z^3)$ for all $z \in [0,1]$ implies that

$$\|x[t] - x_0\|_r = \|y[t] + u(t)h - u(0)h - y[0]\|_r, \text{ for all } t \in (0,T)$$

where $h(z) = \sigma_1 z^2 + \sigma_2 z^3$ for all $z \in [0,1]$. Continuity of the mapping $\Re_+ \ni t \to y[t] \in L_r^2([0,1])$, right continuity of $u : \Re_+ \to \Re$ and the above equality imply that the mapping $[0,T] \ni t \to x[t] \in L_r^2([0,1])$ is continuous. Equations (2.7), (2.8) for $t \in (0,T)$ are verified by using (A.2), (A.8), (A.9), (A.4), (A.5) and (A.7). Moreover, (A.2), (A.3) and (A.8) imply that $x[t] \in C^2([0,1])$ for all $t \in (0,T]$.

Thus we are left with the task of proving that there is a unique function $y: \Re_+ \times [0,1] \to \Re$ is the unique function for which the mapping $\Re_+ \ni t \to y[t] \in L^2_r([0,1])$ is continuous, with $y \in C^1((0,+\infty) \times [0,1])$ satisfying $y[t] \in C^2([0,1])$ for all $t > 0$, $y(0,z) = x_0(z) - u(0)(\sigma_1 z^2 + \sigma_2 z^3)$ for all $z \in [0,1]$ and (A.4), (A.5).

Define:

$$y_0(z) = x_0(z) - u(0)(\sigma_1 z^2 + \sigma_2 z^3) \text{ for all } z \in [0,1] \tag{A.10}$$

$$c_n = \int_0^1 r(z)\phi_n(z)y_0(z)dz, \text{ for } n=1,2,... \tag{A.11}$$

$$\theta_n(t) := \int_0^1 r(z)\phi_n(z)f(t,z)dz, \text{ for all } t \geq 0, \ n=1,2,... \tag{A.12}$$

Since the mapping $[0,1] \ni z \to f(t,z) \in \Re$ is $C^1$ for each $t \geq 0$, it follows from Theorem 11.2.4 in [2], that the following equation holds:

$$f(t,z) = \sum_{n=1}^{\infty} \theta_n(t)\phi_n(z), \text{ for all } (t,z) \in (0,+\infty) \times (0,1) \tag{A.13}$$

Moreover, notice that the Cauchy-Schwarz inequality, in conjunction with the fact that $\|\phi_n\|_r = 1$ (for $n=1,2,...$) and the fact that $f \in C^1(\Re_+ \times [0,1]; \Re)$, implies the following relations:

$$|\theta_n(t)| \leq \left(\int_0^1 r(z)|f(t,z)|^2 dz\right)^{1/2}, \text{ for all } t \geq 0 \tag{A.14}$$

$$\dot\theta_n(t) = \int_0^1 r(z)\phi_n(z)\frac{\partial f}{\partial t}(t,z)dz, \text{ for } t \geq 0 \tag{A.15}$$

$$|\dot\theta_n(t)| \leq \left(\int_0^1 r(z)\left|\frac{\partial f}{\partial t}(t,z)\right|^2 dz\right)^{1/2}, \text{ for } t \geq 0 \tag{A.16}$$

Notice that since the mapping $\Re_+ \ni t \to f(t,z) \in \Re$ is $C^1$, it follows that the mapping $\Re_+ \ni t \to \theta_n(t) \in \Re$ is $C^1$ on $\Re_+$. Since for every $t_1 > 0$ the mapping $[0,t_1] \times [0,1] \ni (t,z) \to f(t,z) \in \Re$ is bounded, there exists $M > 0$ such that $|f(t,z)| \leq M$ for $(t,z) \in [0,t_1] \times [0,1]$. It follows from (A.14) that

$$|\theta_n(t)| \leq M\left(\int_0^1 r(z)dz\right)^{1/2}, \text{ for all } t \in [0,t_1] \text{ and } n=1,2,... \tag{A.17}$$

$$\left|\phi_n(z)\int_0^t \exp(-\lambda_n(t-s))\theta_n(s)ds\right| \leq M\lambda_n^{-1} \max_{0 \leq z \leq 1}(|\phi_n(z)|)\left(\int_0^1 r(z)dz\right)^{1/2},$$
$$\text{for all } (t,z) \in [0,t_1] \times [0,1] \text{ and } n=1,2,... \tag{A.18}$$

Moreover, the Cauchy-Schwarz inequality and (A.11) imply that $|c_n| \leq \left(\int_0^1 r(z)y_0^2(z)dz\right)^{1/2}$ for $n=1,2,...$ and since $\lambda_n \exp(-\lambda_n t) = t^{-1}\lambda_n t \exp(-\lambda_n t) \leq t_0^{-1}\exp(-1)$ for all $t \in [t_0,t_1]$ with $t_0 \in (0,t_1)$, it follows that

$$|\phi_n(z)\exp(-\lambda_n t)c_n| \leq t_0^{-1}\exp(-1)\lambda_n^{-1} \max_{0 \leq z \leq 1}(|\phi_n(z)|)\|x_0\|_r, \text{ for } (t,z) \in [t_0,t_1] \times [0,1], \ n=1,2,... \tag{A.19}$$

Inequalities (A.18), (A.19) and (A.1) imply that the series

$$\sum_{n=1}^{\infty} \phi_n(z)\left(\exp(-\lambda_n t)c_n + \int_0^t \exp(-\lambda_n(t-s))\theta_n(s)ds\right)$$

is uniformly and absolutely convergent on $(t,z) \in [t_0,t_1]\times[0,1]$ for all $t_1 > t_0 > 0$. Therefore, we define $y \in C^0((0,+\infty)\times[0,1])$ by means of the formula:

$$y(t,z) := \sum_{n=1}^{\infty} \phi_n(z)\left(\exp(-\lambda_n t)c_n + \int_0^t \exp(-\lambda_n(t-s))\theta_n(s)ds\right),$$

for all $(t,z) \in (0,+\infty)\times[0,1]$ (A.20)

and we also define

$$y(0,z) := y_0(z), \text{ for all } z \in [0,1] \quad \text{(A.21)}$$

The fact that $y_0 \in L_r^2([0,1])$ (which implies that $\sum_{n=1}^{\infty} c_n^2 < +\infty$) in conjunction with (A.17) the fact that $\sum_{n=1}^{\infty} \lambda_n^{-2}$ (a consequence of (A.1) and the fact that $1 = \int_0^1 r(z)\phi_n^2(z)dz \le \left(\max_{0 \le z \le 1}(|\phi_n(z)|)\right)^2 \int_0^1 r(z)dz$ for $n = 1,2,...$) shows that for all $t \in [0,t_1]$ and for every integer $N \ge 1$, it holds that

$$\|y[t] - y_0\|_r^2 \le 2\sum_{n=1}^{\infty}(\exp(-\lambda_n t)-1)^2 c_n^2 + 2K^2 \sum_{n=1}^{\infty} \lambda_n^{-2}(\exp(-\lambda_n t)-1)^2$$

$$\le 2\sum_{n=N+1}^{\infty} c_n^2 + 2K^2 \sum_{n=N+1}^{\infty} \lambda_n^{-2} + 2(\exp(-\lambda_N t)-1)^2\left(\sum_{n=1}^{N} c_n^2 + K^2 \sum_{n=1}^{\infty} \lambda_n^{-2}\right)$$

where $K := M\left(\int_0^1 r(z)dz\right)^{1/2}$. The above inequality shows that the mapping $\Re_+ \ni t \to y[t] \in L_r^2([0,1])$ is continuous.

From this point the proof is exactly the same with the proof of Theorem 3.1 in [13]. Finally, uniqueness follows from Corollary 2.2 on page 106 of the book [24], since the constructed function $y: \Re_+ \times [0,1] \to \Re$ is a strong solution (see Definition 2.8 on page 109 of the above book). The proof is complete. ◁

**Derivation of (4.31):** Since the mappings $\Re_+ \ni t \to y_n(t)$ for $n = 1,2,...$ are continuous and since differential equations (4.29) hold for all $t \in I$, it follows that the following equations hold for all $t \ge 0$ and $n = 1,2,...$:

$$y_n(t) = \exp(-\mu_n t)y_n(0) - p\frac{d\psi_n}{dz}(1)\int_0^t \exp(-\mu_n(t-s))v(s)ds \quad \text{(A.22)}$$

Equations (A.22) imply that the following equations hold for all $\sigma \in (0,\mu_1)$, $t \ge 0$ and $n = 1,2,...$:

$$|y_n(t)| \le \exp(-\mu_n t)|y_n(0)| + \frac{p}{\mu_n - \sigma}\left|\frac{d\psi_n}{dz}(1)\right|\sup_{0 \le s \le t}\left(|v(s)|\exp(-\sigma(t-s))\right) \quad \text{(A.23)}$$

which, consequently give:

$$|y_n(t)|^2 \le 2\exp(-2\mu_n t)|y_n(0)|^2 + \frac{2p^2}{(\mu_n - \sigma)^2}\left|\frac{d\psi_n}{dz}(1)\right|^2 \sup_{0 \le s \le t}\left(|v(s)|^2 \exp(-2\sigma(t-s))\right) \quad \text{(A.24)}$$

It follows from (4.27), (A.24) and the fact that $0 < \mu_1 < \mu_2 < \ldots < \mu_n < \ldots$, that the following estimate holds for all $\sigma \in (0, \mu_1)$, $t \geq 0$ and $n = 1, 2, \ldots$:

$$\|y[t]\|_2^2 \leq 2\exp(-2\mu_1 t)\|y[0]\|_2^2 + 2p^2 \sup_{0 \leq s \leq t}\left(|v(s)|^2 \exp(-2\sigma(t-s))\right)\sum_{n=1}^{\infty} \frac{1}{(\mu_n - \sigma)^2}\left|\frac{d\psi_n}{dz}(1)\right|^2 \quad (A.25)$$

Changing $z$ by $1-z$ and using Lemma 2.1 in [13], it follows that the boundary value problem

$$p\frac{d^2\bar{x}}{dz^2}(z) - c\bar{x}(z) = 0, \text{ for all } z \in [0,1], \quad (A.26)$$

with

$$b_1\bar{x}(0) + b_2\frac{d\bar{x}}{dz}(0) = 0 \quad , \quad \bar{x}(1) = 1 \quad (A.27)$$

has a unique solution $\bar{x} \in C^2([0,1])$, which satisfies

$$p^2 \sum_{n=1}^{\infty} \mu_n^{-2}\left|\frac{d\psi_n}{dz}(1)\right|^2 = \int_0^1 \bar{x}^2(z)dz \quad (A.28)$$

Consequently, we obtain from (A.28):

$$K := p^2 \sum_{n=1}^{\infty} (\mu_n - \sigma)^{-2}\left|\frac{d\psi_n}{dz}(1)\right|^2 \leq \left(\frac{\mu_1}{\mu_1 - \sigma}\right)^2 \int_0^1 r(z)\bar{x}^2(z)dz < +\infty \quad (A.29)$$

Therefore, we obtain from (A.25) and (A.29) that the following estimate holds for all $\sigma \in (0, \mu_1)$ and $t \geq 0$:

$$\|y[t]\|_2^2 \leq 2\exp(-2\sigma t)\|y[0]\|_2^2 + 2K \sup_{0 \leq s \leq t}\left(|v(s)|^2 \exp(-2\sigma(t-s))\right) \quad (A.30)$$

Inequality (4.31) is a direct consequence of inequality (A.30). The derivation is complete. ◁